\documentclass[11pt]{amsart}
\usepackage{amsfonts,amssymb}
\usepackage{graphicx}
\usepackage{amscd}
\usepackage{eucal}
\usepackage{amsmath}
\usepackage{eufrak}
\input epsf
\newtheorem{thm}{\bf Theorem}[section]
\newtheorem{cor}[thm]{\bf Corollary}

\newtheorem{prop}[thm]{\bf Proposition}

\def\LL{\rm{L}}
\def\R{\mathbb{R}}

\def\Q{\mathcal{Q}}
\def\S{{\mathcal S}}
\def\H{{\mathcal H}}

\def\SO{{\rm SO}}

\def\<{\langle }
\def\>{\rangle }
\def\a{{\mathbf a}}
\def\v{{\mathbf v}}

\def\L{\mathcal{L}}
\def\D{\mathcal{D}}

\def\u{\mathbf{u}}
\def\v{\mathbf{v}}
\def\x{\mathbf{x}}
\def\z{\mathbf{z}}
\def\y{\mathbf{y}}
\def\w{\mathbf{w}}

\def\at{\textsf{a}}
 \def\bt{\textsf{b}}
\def\ct{\textsf{c}}
\def\mt{\textsf{m}}
\def\nt{\textsf{n}}
\def\rt{\textsf{r}}
\def\ot{\textsf{o}}
\def\hh{\textsf{h}}
\def\doth{\ \bar{\cdot}\ }
\def\crossh{\  \bar{\wedge} \ }

\begin{document}
\title{Projective geometry  from  Poisson  Algebras}
\author{Francesca Aicardi}

\begin{abstract}
In  analogy  with  the  Poisson  algebra  of  the quadratic  forms  on  the
symplectic  plane,  and the notion of  duality  in the  projective  plane   introduced  by  Arnold  in  \cite{Arn},  where
the  concurrence  of  the  triangle  altitudes  is  deduced  from
the  Jacobi  identity, we  consider  the  Poisson  algebras  of  the  first  degree  harmonics
on  the  sphere,  the pseudo-sphere  and  on the
hyperboloid,  to  obtain  analogous  duality  notions  and
similar results  for  the  spherical,  pseudo-spherical  and  hyperbolic
geometry.    Such  algebras, including  the  algebra  of  quadratic  forms,  are    isomorphic, as  Lie
algebras, either to  the  Lie algebra  of  the  vectors  in  $\R^3$,  with  vector  product, or to algebra $sl_2(\R)$.   The  Tomihisa
identity, introduced in  \cite{Tom} for the algebra  of  quadratic  forms,  holds   for all
these    Poisson  algebras and   has  a  geometrical
interpretation.  The  relation  between  the different definitions  of
 duality  in   projective  geometry  inherited  by  these
 structures  is   shown.
\end{abstract}
\maketitle

\section*{Introduction}

In \cite{Arn},   Arnold  has  shown that a hyperbolic  version of the concurrence   altitudes
theorem  for  triangles  holds
in   the projective  plane,   to which the  space  of  binary  quadratic forms  projects.
In  fact,   this  theorem is  a direct consequence  of  the  Jacobi  identity  for  the
Poisson  brackets, in the Poisson  algebra of the binary quadratic  forms $ap^2+2bpq+cq^2$   on the
symplectic plane  $(p,q)$.

The space  of the  coefficients $(a,b,c)$  of the quadratic forms  is also endowed in  \cite{Arn}
with  a scalar product (defining its metrics)  as  well as  with  a vector product  (the Poisson  bracket),  so  that   the orthogonality  of two  forms is defined as  the vanishing  of  their scalar product,  and the
Poisson bracket of two forms  is  a third  form orthogonal  to the  two initial forms.

 A  geometrical notion  of  duality   between  points and  lines in  the  projective plane
  is   also  introduced in  \cite{Arn}, which  allows
to give  a  geometrical  meaning  to  any  expression  involving  Poisson  brackets.

The  starting point  of  the present work was the  following  consideration.

The  real vector  space  $\R^3$  is  a Lie   algebra  if  endowed  with  the vector  product
\begin{equation}\label{vectorp} (x,y,z)\wedge(x',y',z')=(yz'-zy',zx'-xz',xy'-yx').  \end{equation}
The  vector product (\ref{vectorp}) is  $\SO(3)$ invariant  and  has  indeed   a  geometrical  meaning
in the oriented $\R^3$,  provided with  Euclidean metric.

However,  the  six  vector products  in $\R^3$,  corresponding  to    choices of the signs $\sigma_i$ for the components of   (\ref{vectorp})  such  that  not  all  signs  are  coinciding
\begin{equation}\label{VPS} (x,y,z)\wedge(x',y',z')=\left( \sigma_1(yz'-zy'),\sigma_2(zx'-xz'),\sigma_3(xy'-yx')\right) \end{equation}
satisfy  the Jacobi identity  and  the corresponding Lie  algebras are all isomorphic  to the algebra  $sl_2(\R)$.

Each of  these vector products  has a  geometrical meaning  in   $\R^3$,  endowed  with an orientation
and  with the (pseudo)metric $g_{i,j}=\sigma_i\delta_{i,j} $.  The  vector  product of two vectors $\v$ and $\v'$ is orthogonal  to both $\v$ and $\v'$,
the orthogonality being the  vanishing of the scalar  product defined  by  the metric.

Now,  the Poisson algebra of quadratic forms, together  with  the scalar  product introduced  by  Arnold,
is,  up  to  a change  of  coordinates,   the Lie  algebra  of $\R^3$  with   vector product  (\ref{VPS})
and  pseudo Euclidean  metric
given by $\sigma_1=\sigma_2=1,\sigma_3=-1$.

A  natural  question is: does it  exist a Poisson algebra  on a  symplectic  manifold corresponding  to  the
Euclidean  $\R^3$?

In  Section 1  we show that  such  algebra is the  3-space  of  spherical harmonics  of  first degree,  with
$\LL^2$-metric.

Projecting  this  space to the  unit 2-sphere, we  associate  with each  harmonic  a pair  of antipodal  points
and  its  {\it dual} object,  the great circle equidistant  from  these  points.  Like in  \cite{Arn}, the
Poisson bracket  has   a  geometrical interpretation  and the  Jacobi identity
assumes the meaning of the altitudes concurrence theorem  for  spherical  triangles.

In  Section 2 we  show that the  pseudo Euclidean  geometry of the  Poisson  algebra of quadratic forms can  be
obtained  as  well as  the Poisson  algebras   of the  hyperbolic harmonics  of first degree   on  the
one sheet  and on the two sheeted  hyperboloids.   In these  cases  the  scalar product  and  the metrics
are not  defined   intrinsically  (as   in  the  spherical  case,  by the $\LL^2$ norm), but  the  vanishing
of  the  scalar  product of two harmonics  has    a  geometrical  meaning  in terms of   the  metrics
  of the  hyperboloids  similar to the meaning of the  vanishing
of  the  scalar  product of two spherical harmonics  in  terms  of  the metric  of the  sphere.

In Section 3  we  see  the relation  between  the projective geometries obtained  projectivizing    the
pseudo-Euclidean  space  (as in  \cite{Arn}) and  the Euclidean  one.  Moreover, we  observe  that
in the  case of pseudo-Euclidean  geometry,  any non degenerate conic  in  the plane can  play the  role of
the absolute of the Lobachevsky disc  for  a  suitable projection: the orthogonality  and  duality  can  be therefore
defined using this conic, obtaining   the  concurrence  theorem for  the  corresponding triangle `altitudes', defined
in this  projection.

Following  the Arnold suggestion:  ``{\it One might  use the Jacobi identity and other  theorems  of the
quadratic  forms symplectic algebra to obtain new results of projective geometry}'',   Tomihisa \cite{Tom}
has recently shown  one identity  holding  for  the Poisson brackets in  this  algebra,
from which  all basic theorems  of Projective Geometry  follow.

In  Section  4  we  show  that  the Pappus  theorem  itself  is
contained  in  the  Tomihisa  identity.  We observe,  moreover,   that  the  geometrical
interpretation  of  the  Tomihisa identity
in projective  geometry is  independent from the  metric  structure  and holds  in  the  Lie-algebra
of the  Euclidean   $\R^3$  as well.

\section*{Notations  and  preliminary  facts}

In this  paper  we  use the following notations:

$\R^3$ is the  oriented  three-dimensional vector  space.

$\R^3_E$  is  the Lie  algebra of  the vectors  in   the Euclidean  3-space, with  vector product (\ref{vectorp}).
 The   scalar  product of   two  vectors  $\v=(x,y,z)$ and $\v'=(x',y',z')$ is:
\begin{equation}\label{scalarp}  \v \cdot \v'=xx'+yy'+zz',
  \end{equation}
and  the square norm of  $\v$  is
\begin{equation}\label{norm}  ||\v||:=\v\cdot \v= x^2+y^2+z^2. \end{equation}

$\R^3_H$  is  the  Lie algebra of the  vectors in the  the pseudo-Euclidean  3-space, with metric   $g_{i,j}=\sigma_i \delta_{i,j}$,  $\sigma_1=\sigma_2=1,\sigma_3=-1$.  I.e.,  the   vector product  of   two  vectors  $\v=(x,y,z)$ and $\v'=(x',y',z')$ is:
\begin{equation}\label{vectorph}  \v \bar \wedge \v' =\left( (yz'-zy'),(zx'-xz'),-(xy'-yx')\right), \end{equation}
their  {\it scalar  product}
\begin{equation}\label{scalarph}  \v \ \bar{\cdot}\   \v'=xx'+yy'-zz',
  \end{equation}
and the  {\it square  norm}  of    $\v$  is
\begin{equation}\label{normh} \bar{||}\v \bar{||}:=\v  \doth \v= x^2+y^2-z^2. \end{equation}

\noindent{\it Remark 1.} Observe that   $\v \wedge \v' \cdot \v= \v \wedge \v' \cdot \v'=0$ and $\v \crossh \v' \doth \v= \v \crossh \v' \doth \v'=0$. Moreover,  $\v \wedge \v'=0$ ($\v \crossh \v'=0$) iff  $\v'=\lambda \v$  for a real $\lambda$.

\noindent{\it Remark 2.} Given  three  independent  vectors in  the oriented Euclidean  3-space,  the  scalar product  of any one of them (say, $\v$) with the  vector  product  of the two  others  vectors  (say, $\u,\w$) is  equal  to  the  volume  of  the parallelepiped  defined by  the three  vectors,  with
positive (negative)  sign  if  the  ordered triple $(\v,\u,\w)$  orients $\R^3$  positively (negatively).
The  same  holds  in $\R^3_H$.  This allows to give  a  definition of  the  vector  product  in $\R^3_H$,
independent  of the  coordinates and based only on the  $\SO(2,1)$-invariant metric:  the  vector product
of  two  vectors  $(\u,\w)$  is  the  vector  such  that  its  scalar  product with  any third  vector  $\v$
lying  outside the plane of $\u$ and  $\w$  is  equal  to the  oriented  volume  defined  by  $(\u,\w,\v)$.

\noindent{\it Remark 3.}   The algebra  $R^3_H$   is isomorphic to  the algebra of  trace  zero  $2\times 2$ real  matrices.
 The isomorphism  $\mu:R^3_H\rightarrow sl_2(\R)$ is   given  by
 \[    \mu (X_1,X_2,X_3) = \frac{1}{2} \begin{pmatrix} X_2  &  -(X_1+X_3) \\
X_3-X_1  &  -X_2  \end{pmatrix}.  \]

We  will consider  also:

\noindent {\it i.}
$S^2\subset \R^3_E$,  the    sphere  $x^2+y^2+z^2=1$, with coordinates
$\theta$ and  $\phi$,  related to  $x,y,z$ by:
\[    x:=\sin \theta \cos \phi, \quad  y:=\sin \theta \sin \phi, \quad  z:=\cos \theta.  \]

The  restriction of the  metric of $\R^3_E$  to the   sphere   yields  to the metric
\begin{equation}\label{metricS} ds^2= d\theta^2+ \sin^2 \theta d\phi^2. \end{equation}

\noindent {\it ii.}
$L^2\subset \R^3_H$,  the    pseudo-sphere  $x^2+x^2-z^2=-1$, $z>0$, with coordinates
$\chi$ and  $\phi$,  related to  $x,y,z$ by:
\[    x:=\sinh \chi \cos \phi, \quad  y:=\sinh \chi \sin \phi, \quad  z:=\cosh \chi.  \]

The  restriction of the  pseudometric of $\R^3_H$  to the hyperboloid $L^2$    yields to the metric

\begin{equation}\label{metricL} ds^2= d\chi^2+ \sinh^2 \chi d\phi^2.  \end{equation}

\noindent {\it iii.}
$D^2\subset \R^3_H$,  the one-sheeted hyperboloid   $x^2+y^2-z^2=1$, with coordinates
$\chi$ and  $\phi$,  related to  $x,y,z$ by:
\[    x:=\cosh \chi \cos \phi, \quad  y:=\cosh \chi \sin \phi, \quad  z:=\sinh \chi.  \]

The  restriction of the  pseudometric of $\R^3_H$  to the hyperboloid  $D^2$  yields  to the metric
  \begin{equation}\label{metricD} ds^2= -d\chi^2+ \cosh^2 \chi d\phi^2.  \end{equation}
Note that  metric  (\ref{metricD}) is indefinite.

\noindent{\it  Remark 4.} A  coordinates-independent  description  of the  metrics   (\ref{metricL})
and (\ref{metricD})  is given  in  \cite{Arn2}.

\section{Spherical geometry and  the algebra $\S_1$}

In this  section  we obtain the  spherical  geometry  from  the Poisson algebra  of  the  spherical  harmonics
of degree 1.

Consider  the space $\mathcal S_1$  of the linear combinations with  real  coefficients of the following  spherical  harmonics of  first  degree,
\[  f_1:=   \sin \theta  \cos \phi;  \quad    f_2:=  \sin \theta \sin \phi; \quad  f_3:=  \cos \theta, \]
\[  \mathcal S_1=\{ Xf_1+Yf_2+Zf_3,  \quad (X,Y,Z)\in \R^3   \}. \]

 The $\LL^2$-metric with respect to the standard area form on the sphere  defines the scalar product:

\begin{equation}\label {prodL2}\langle F, F'\rangle=  \frac{3}{4\pi} \int_{S^2}  FF'  \sin\theta d\theta d \phi.  \end{equation}

We identify   the  function $F=Xf_1+Yf_2+Zf_3$  with  the  triple of its coefficients:  $F=(X,Y,Z)\in \R^3$.

\begin{prop} \label{prop1.1}  The scalar product (\ref{prodL2}) coincides   with  the  scalar product  (\ref{scalarp}) in  $\R^3_E$.
\end{prop}

\noindent{\it Proof.}  The three  functions:
 $ f_1$, $f_2$ and   $f_3$
 constitute  a orthonormal  basis:
 \[ \langle f_i , f_j \rangle =  \delta _{ij}.  \]
\hfill $\square$

\noindent{\it  Definition.}  Two functions in $\S_1$  are {\it orthogonal}  iff their scalar  product  vanishes.

Consider the Poisson structure associated with the
standard symplectic form on the sphere:
\[  \{F,F'\}=  \frac{1}{\sin\theta} \left(\frac{\partial{F}}{\partial \theta}  \frac{\partial{F'}}{\partial \phi}-  \frac{\partial{F}}{\partial \phi} \frac{\partial{F'}}{\partial \theta}\right).  \]

\begin{prop} \label{prop1.2} The   Poisson bracket of two   functions   in  $\S_1$  satisfies:
\[ \{F,F'\}=  F \wedge F'.  \]
\end{prop}
\noindent{\it  Proof.}  If  $F=(X,Y,Z)$ and  $F'=(X',Y',Z')$ are in  $\S_1$, then $\{F,F'\}\in \S_1$.  For
basic functions  $f_i$,  $i=1,2,3$,  we  get
\[  \{f_1,f_2 \}=- \{f_2,f_1\}=f_3, \quad    \{f_2,f_3\}= -\{f_3,f_2\}=f_1,\quad    \{f_3,f_1\}= -\{f_1,f_3\}=f_2.  \]
Since  the  Poisson  bracket  is  bilinear  and antisymmetric, Proposition \ref{prop1.2} follows.  \hfill $\square$

\noindent{\it Definition.} Two functions  having
non zero  Poisson  bracket are  said  {\it independent}.

We  see  now  a  geometrical implication of the  orthogonality  of  two
spherical functions of $ \S_1$.

\begin{prop} \label{prop1.3}  Let  $F_1=X_1 f_1+Y_1f_2+ Z_1f_3$  and  $F_2=X_2 f_1+Y_2f_2+ Z_2f_3$ be  two non  zero
spherical functions. If  $F_1 \cdot F_2=0$,   then  the great circles
 $F_1=0$  and $F_2=0$  on the  sphere  meet  orthogonally.
\end{prop}

\noindent{\it  Proof.}  Let  $s$  be the parameter  along the curves on the  sphere where  $F_1$ and  $F_2$ vanish:
\[ F_1(\theta_1(s),\phi_1(s))=0,  \quad  F_2(\theta_2(s),\phi_2(s))=0,  \]
so that  at  $s=s^*$, $\theta_1=\theta_2:=\theta^*$  and   $\phi_1=\phi_2:=\phi^*$.
Since (\ref{metricS})  holds, we have  to  prove that  at  $s=s^*$   the  following equation  is fulfilled, whenever $X_1X_2+Y_1Y_2+Z_1Z_2=0$:
\begin{equation}\label{met}    \frac{d\theta_1}{ds}\big|_{s=s^*} \frac{d\theta_2}{ds}\big|_{s=s^*}+  \sin^2\theta^* \frac{d\phi_1}{ds}\big|_{s=s^*}\frac{d\phi_2}{ds}\big|_{s=s^*}=0.
\end{equation}
Take  as  parameter $s$ the  angle $\phi$.
The functions $F_i$  vanish on  the curves
\[ \theta_i=\arctan\left(\frac {Z_i}{X_i \cos\phi+Y_i\sin\phi}\right). \]
From the  equation   $\theta_1(\phi^*)=\theta_2(\phi^*)$ we obtain
\[  \phi^*= \pm\arctan(({X_1Z_2-X_2Z_1})/({Y_1Z_2-Y_2Z_1}),\]
and
 \[  \theta^*=\arctan\left(\frac{\sqrt{(X_1Z_2-X_2Z_1)^2+(Y_1Z_2-Y_2Z_1)^2}}{(Y_1X_2-X_1Y_2)}\right). \]

The left member of Eq. \ref{met} becomes:
\begin{equation}    \frac{d\theta_1}{d\phi}\big|_{\phi=\phi^*} \frac{d\theta_2}{d\phi}\big|_{\phi=\phi^*}+  \sin^2 \theta^*= \frac{(X_1Z_2-X_2Z_1)^2+(Y_1Z_2-Y_2Z_1)^2}{Z_1Z_2 \  || \{ F_1, F_2 \}||}(F_1\cdot F_2),    \end{equation}
which  vanishes  if  $F_1\cdot F_2=0$. \hfill $\square$

\subsection{From  Algebra to Geometry}\label{sphere}
We  define  the following duality on  the  unit  2-dimensional sphere.

\noindent{\it  Definition.}   A pair of  antipodal  points  is {\it dual}  of the   great circle
equidistant  from  these points.

For example,  the  pair of  North-South poles  is  dual  of the equator.

We  denote  by  Greek letters  (e.g. $\alpha$)   the  great  circles,  by  small letters (e.g. $a$) the
pairs  of  antipodal  points,  and by       $(  a  |  \alpha  )$  the  pair  of  dual  objects  $a$  and  $\alpha$.

With every  vector ${\a}=(X,Y,Z)$  in  $\R^3$  we  associate the  pair  of antipodal points
  \[  a= (  \at, -\at ),  \quad   \at= \frac{\a}{\sqrt{||\a ||}} \]
  on the unit  sphere.

Every function  $F\in \S_1$  is  represented  by  a  vector  $\a=(X,Y,Z)$, therefore  we
associate with every function  $F\in \S_1$  the  corresponding pair  $ (  a  |  \alpha  ) $.

\begin{prop} \label{prop1.4}  Let  $F,F'\in  \S^1$ be two independent  functions, and let   $\{F,F'\}=F''$.
Let  $\a,\a',\a''$  be the  corresponding  vectors, $a,a',a''$   the  corresponding pairs  of antipodal points on the  sphere  and $\alpha,\alpha',\alpha''$  their  dual  great circles.   Then
 $  (  a''  |  \alpha''  ) $
has the  following geometrical meaning:   the  antipodal points of the pair  $a''$ are the intersection of  the  two  great  circles
$\alpha$ and  $\alpha'$.  The circle   $\alpha''$  is the great circle  joining  the two  pairs of points in  $a$  and  $a'$.   \end{prop}

\noindent{\it Proof.}  The  vector  $\a''$ is, by Proposition \ref{prop1.2}, the vector product of $\a$ and  $\a'$, which is
orthogonal to the plane containing     $\a$ and  $\a'$.  The great circle containing $a$  and  $a'$
is therefore  the  intersection  of this plane with the sphere, which is the circle  $\alpha''$  dual of $a''$.
The  circles  $\alpha$  and  $\alpha'$  lie  in two  planes  orthogonal respectively to  the  vectors  $\a$  and  $\a'$.  Since $\a''$ is orthogonal  to both  $\a$  and  $\a'$,  $\a''$
lies in the  intersection of these planes, and therefore the pair $a''$  is the intersection of $\alpha$ with  $\alpha'$. \hfill  $\square$

\noindent{\it  Definition.} A  spherical triangle is  said  to be  {\it proper},  if  no one of its  vertices  belong to  the pair of points dual of the  great  circle  containing the other vertices,  or,  equivalently,
no  one of its  sides belongs to a  great  circle dual  of  a pair containing the opposite vertex.

 \begin{prop} \label{prop1.5} A  spherical  triangle  defined  by the  normalized   vectors  $\at_1$, $\at_2$, $\at_3$, corresponding  to  the  functions  $F_1,F_2,F_3 \in \S_1$
is  proper  iff  $\{\{F_i,F_j\},F_k\}\not=0$  for  any  choice of the indices  $i,j,k$  among  the  permutations
of $1,2,3$.\end{prop}

\noindent{\it Proof.} The  identity  $\{\{F_i,F_j\},F_k\}=0$  holds if and only if  $\{F_i,F_j\}=\lambda F_k$,  by Remark 2.
If  $\{F_i,F_j\}=\lambda F_k$  for some  $i,j,k$,  then the  normalized  vector $\at_3$ is orthogonal to the plane of  $\at_1$ and  $\at_2$, i.e., $a_3$  is  dual  of  the  great circle  containing  the pairs $a_1$ and  $a_2$. \hfill $\square$

\noindent{\it  Definition.}  The {\it  altitude}  from a  vertex   $\at_i$ of a  proper spherical  triangle  is the  great
circle  through $\at_i$  intersecting  orthogonally   the great  circle  containing  the  opposite  side to  $\at_i$.

\begin{prop} \label{prop1.6} The  altitude from a  vertex   $\at_i$ of a  proper spherical  triangle  is the
great circle  through the pair  $a_i$   passing trough the pair  $b$,  dual  of   the  the great
circle $\beta$ containing  the  other two vertices of the  triangle. \end{prop}

\noindent{\it Proof.}  Any circle  passing  through $b$  is orthogonal  to  the great  circle $\beta$, by the
definition of  duality. \hfill $\square$

\noindent{\it  Remark 5.}  The  altitude  from a  vertex  $\at_i$  of  a proper spherical triangle is  also
the altitude of  the antipodal  triangle (with vertices the antipodal  $-\at_1,-\at_2,-\at_3$),
and  of the  other six triangles cut off  on the  sphere by  the  great circles,  sides of  the
triangle.

\begin{prop} \label{prop1.7}  Let  $F_1,F_2,F_3 \in \S_1$ define   the  vertices  of  a  proper  spherical  triangle.
Let  $\at_1,\at_2,\at_3$  be the normalized vectors  corresponding  to $F_1,F_2,F_3$, $a_1,a_2,a_3$
 the  corresponding pairs  of antipodal points on the  sphere  and $\alpha_1,\alpha_2,\alpha_3$   their  dual
 great circles.
Then  the  Jacobi identity:
 \begin{equation}\label{Jak}  \{\{F_1,F_2\},F_3\}+  \{\{F_2,F_3\},F_1\}+  \{\{F_3,F_1\},F_2\}=0 \end{equation}
has the  following geometrical meaning:
i) The altitudes   of  the  eight  spherical triangles  with vertices  $\pm \at_1,\pm \at_2, \pm \at_3$
meet  at the same  pair $h:=(\hh,-\hh)$ of  points;
ii)  These altitudes  are  also  the  altitudes of the  eight  spherical triangles defined  by  the three
great circles  $\alpha_1,\alpha_2,\alpha_3$. \end{prop}

\noindent{\it  Proof.} Equation  (\ref{Jak})  says that  the  three  vectors $\{\{F_1,F_2\},F_3\}$,
 $\{\{F_2,F_3\},F_1\}$ and   $\{\{F_3,F_1\},F_2\}$   have zero sum.   In particular, they lie in the same plane.
{\it i}) By Proposition  \ref{prop1.4},   $\{\{F_1,F_2\},F_3\}$  represents the  great circle
$\gamma_3$ containing $\at_3$ and  the  pair  $b_3$ of points, dual of  the opposite  side $\beta_3$, which is  the
great  circle   containing  $\at_1$  and  $\at_2$.  But $\{\{F_1,F_2\},F_3\}$   represents also the  pair of
points  $c_3$ (dual of  $\gamma_3$), which is the  intersection  of  $\alpha_3$  with       $\beta_3$.    The equation
implies that  the pair $c_3$  and the  analogue  $c_1$  and  $c_2$  lie  in the  same great  circle. Therefore,
the  dual  great circles  $\gamma_3$, $\gamma_1$  and  $\gamma_2$  meet on a  unique pair of  points  $h$.
The  three altitudes $\gamma_1$, $\gamma_2$ and  $\gamma_3$
are common to all  eight spherical triangles,  cut off  by  the  3 sides  of the  initial  triangle.

\begin{figure}[h]
\centerline{\epsfbox{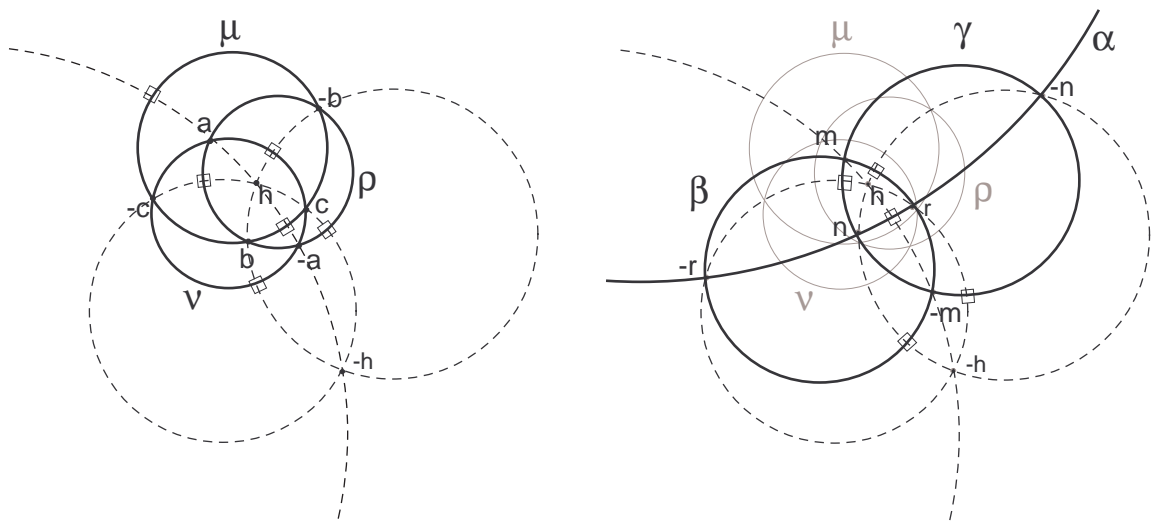}} \caption{Stereographic projection  of  the spherical
triangles defined  by the great circles  $\mu, \nu,\rho$  through  $\at,\bt,\ct$,  and  the  dual
triangles  with great circles  $\alpha,\beta,\gamma$  trough $\mt,\nt,\rt$.   They have  the  same
altitudes (dotted lines) meeting  at $\pm \hh$.    }\label{proj03}
\end{figure}

{\it ii})  The  spherical  triangles   defined  by the three  dual   great  circles  $\alpha_1,\alpha_2,\alpha_3$
have the  vertices  at the  points  of the pairs $b_1,b_2,b_3$,   dual  of  the sides of  the  initial triangle.
The  altitudes  are  the great  circles  containing  each one  a pair  $b_i$  and  the  corresponding  pair $a_i$, dual of  $\alpha_i$.  Therefore  these  altitudes coincide  with  the great circle  $\gamma_i$, meeting at the pair of points  $h$ (See Figure \ref{proj03}). \hfill  $\square$

\noindent{\it Remark 6.} Observe  that, if  the  spherical  triangle $\at_1,\at_2,\at_3$ with sides  $\beta_1, \beta_2,\beta_3$ is  not  proper, that  is  a  vertex,  say $\at_1$, belongs  to the  pair of  points  dual  to the  great circle $\beta_1$  through  $\at_2$  and  $\at_3$ (i.e.,  $\beta_1=\alpha_1$)
then  all  great  circles  through $\at_1$  are altitudes  of the  triangle,   and  both   the  altitude $\gamma_2$  from  $\at_2$ to  $\beta_2$  and the  altitude  $\gamma_3$  from  $\at_3$  to $\beta_3$ coincide  with the  side  $\alpha_1$.
The  common point  of  the  altitudes of  the  triangle  degenerates  in this  case  to  the
entire  circle $\alpha_1$.   This  is  the meaning of the   Jacobi  identity  that in this case
reads   $ \{\{F_1,F_2\},F_3\}+     \{\{F_3,F_1\},F_2\}=0$.  Indeed,   the  points  $ \{\{F_1,F_2\},F_3\}$  and    $ \{\{F_3,F_1\},F_2\}$, being antipodal,  define two coinciding pairs of
points, $c_2$ and  $c_3$,  dual to the  altitudes $\gamma_2$ and  $\gamma_3$, that therefore coincide.

\section{Hyperbolic  geometry}

\subsection{Lobachevsky geometry}

In  analogy with  the  Poisson algebra  $\S_1$, we consider  the Poisson algebra $\mathcal H^-_1$  of the linear combinations with  real  coefficients of the following   hyperbolic harmonics of  first  degree,  defined on the upper sheet  of the two sheeted hyperboloid  (the psedo-sphere):
\[  f_1:=   \sinh \chi  \cos \phi;  \quad    f_2:=  \sinh \chi \sin \phi; \quad  f_3:=  \cosh \chi, \]
\[  \mathcal H^-_1=\{ Xf_1+Yf_2+Zf_3,  \quad (X,Y,Z)\in \R^3   \}. \]

The Poisson bracket in $\H^-_1$  is:
\begin{equation}\label{PoissonL} \{F,F'\}=  \frac{1}{\sinh\chi} \left(\frac{\partial{F}}{\partial \chi}  \frac{\partial{F'}}{\partial \phi}-  \frac{\partial{F}}{\partial \phi} \frac{\partial{F'}}{\partial \chi}\right).  \end{equation}

\begin{prop} \label{prop2.1}  The Poisson  algebra $\H^-_1$   coincides with the Lie
algebra  $\R^3_H$.  \end{prop}
\noindent{\it  Proof.} If $F=(X,Y,Z)$ and  $F'=(X',Y',Z')$ are in  $\H^-_1$, then $\{F,F'\}\in \H^-_1$,  and coincides with the  vector product (\ref{vectorph})  in $\R^3_H$:
\begin{equation} \label{eqPL}  \{F,F'\}=  F \crossh F'.  \end{equation}
Indeed, consider the
basic functions  $f_i$,  $i=1,2,3$, and  their  Poisson  brackets.  They satisfy:
\[  \{f_1,f_2 \}=- \{f_2,f_1\}=-f_3, \quad    \{f_2,f_3\}= -\{f_3,f_2\}=f_1,\quad    \{f_3,f_1\}= -\{f_1,f_3\}=f_2.  \]
Since  the  Poisson  bracket  is  bilinear  and antisymmetric,  Eq.  (\ref{eqPL}) is  fulfilled. Proposition \ref{prop2.1} follows.  \hfill $\square$

Also in  $\H^-_1$ two  functions  having
non zero  Poisson  bracket will  be   said  {\it independent}.

The Poisson  algebra  $\H^-_1$ can  be  therefore  endowed  with  the same metric  as $\R^3_H$, with scalar product
(\ref{scalarph}),   square norm  (\ref{normh}),  and  the  same  definition of  orthogonality.
    In this  way  the  Poisson  bracket  of  two independent  functions
is orthogonal  to  these functions.

\noindent{\it  Remark 7.}  The   three functions  $f_1$,  $f_2$  and  $f_3$ are orthogonal  and  satisfy   $f_1^2+f_2^2-f_3^2=-1$. They are the   coordinates in  $\R^3$ of   the point with  coordinates $(\chi, \phi)$  on the  pseudo-sphere.

Also  in this  case we see the  geometrical implication of the  orthogonality  of  two
functions of $\H^-_1$.

\begin{prop} \label{prop2.2}  Let  $F_1=X_1 f_1+Y_1f_2+ Z_1f_3$  and  $F_2=X_2 f_1+Y_2f_2+ Z_2f_3$ be  two non  zero
functions of  $\H_1^-$. If  $F_1 \doth F_2=0$,   then  the lines
 $F_1=0$  and $F_2=0$  on the  pseudo-sphere  meet  orthogonally.
\end{prop}

\noindent{\it  Proof.}  The  proof is  analogous of that of  \ref{prop1.3}.
If $\chi_1(\phi)$ and  $\chi_2(\phi)$  parametrize  the  curves $F_1=0$  and  $F_2=0$  which meet  at the point
 $(\chi^*, \phi^*)$, we    obtain finally,  using  metric (\ref{metricL}):
\begin{equation}    \frac{d\chi_1}{d\phi}\big|_{\phi=\phi^*} \frac{d\chi_2}{d\phi}\big|_{\phi=\phi^*}+  \sinh^2 \chi^*=\frac{(X_1Z_2-X_2Z_1)^2+(Y_1Z_2-Y_2Z_1)^2}{Z_1Z_2\  \bar{||} \{ F_1, F_2 \}\bar{||}}(F_1\doth F_2),    \end{equation}
which  vanishes  if  $F_1\doth F_2=0$. \hfill $\square$

\subsection{De Sitter geometry}

In complete  analogy   with  the algebra $\mathcal  H^-_1$, we consider   the Poisson  algebra $\mathcal H^+_1$  of the linear combinations with  real  coefficients of the following   hyperbolic harmonics of  first  degree,  defined on the one sheeted   hyperboloid:
\[  f_1:=   \cosh \chi  \cos \phi;  \quad    f_2:=  \cosh \chi \sin \phi; \quad  f_3:=  \sinh \chi, \]
\[  \mathcal \H^+_1=\{ Xf_1+Yf_2+Zf_3,  \quad (X,Y,Z)\in \R^3   \}, \]
with Poisson bracket
\begin{equation}\label{PoissonD}  \{F,F'\}=  \frac{1}{\cosh\chi} \left(\frac{\partial{F}}{\partial \chi}  \frac{\partial{F'}}{\partial \phi}-  \frac{\partial{F}}{\partial \phi} \frac{\partial{F'}}{\partial \chi}\right).
\end{equation}

\begin{prop} \label{prop2.3} The  Poisson   algebra $\H^+_1$ is  isomorphic to  the  Lie algebra   $\R^3_H$;  the  isomorphism
is  given by an  inversion of  sign: $F \rightarrow -F$.
\end{prop}
\noindent{\it  Proof.}    We prove that if  $F=(X,Y,Z)$ and  $F'=(X',Y',Z')$ are in  $\H^+_1$, then $\{F,F'\}\in \H^+_1$,  and
is opposite to  the  vector  product  (\ref{vectorph}) in  $\R^3_H$.
\begin{equation}\label{eqPD} \{F,F'\}= - F\crossh F' .  \end{equation}   Consider the
basic functions  $f_i$,  $i=1,2,3$, and  their  Poisson  brackets.  We  obtain
\[  \{f_1,f_2 \}=- \{f_2,f_1\}=f_3, \quad    \{f_2,f_3\}= -\{f_3,f_2\}=-f_1,\quad    \{f_3,f_1\}= -\{f_1,f_3\}=-f_2.  \]
Since  the  Poisson  bracket  is  bilinear  and antisymmetric, Eq. (\ref{eqPD})  is fulfilled.  Therefore  $\{-F,-F'\}=-F\crossh F'$ and the Proposition  follows.  \hfill $\square$

Also  here the  same metric    structure  as in  $R^3_H$  can be introduced,  so  that  the Poisson bracket  of
  two  functions in $\H^+_1$ is orthogonal to  both functions.

\noindent{\it  Remark 8.}  The three functions  $f_1$,  $f_2$  and  $f_3$ are orthogonal and satisfy   $f_1^2+f_2^2-f_3^2=1$. They are
the   coordinates in  $\R^3$ of   the point with  coordinates $(\chi, \phi)$  on the one sheeted hyperboloid.

The  geometrical meaning of the orthogonality of two functions in  $\H^+_1$ is the following.

\begin{prop} \label{prop2.5}  Let  $F_1=X_1 f_1+Y_1f_2+ Z_1f_3$  and  $F_2=X_2 f_1+Y_2f_2+ Z_2f_3$ be  two non  zero
functions in  $\H_1^+$. If  $F_1 \cdot F_2=0$,   then  the lines
 $F_1=0$  and $F_2=0$  on the  hyperboloid  meet  orthogonally.
\end{prop}

\noindent{\it  Proof.}  The  proof is  analogous of that of  Proposition \ref{prop1.3}.
If $\chi_1(\phi)$ and  $\chi_2(\phi)$  parametrize  the  curves $F_1=0$  and  $F_2=0$  which meet  at the point
 $(\chi^*, \phi^*)$, we    obtain finally, using  metric (\ref{metricD}):
\begin{equation}    -\frac{d\chi_1}{d\phi}\big|_{\phi=\phi^*} \frac{d\chi_2}{d\phi}\big|_{\phi=\phi^*}+  \cosh^2 \chi^*= - \frac{(X_1Z_2-X_2Z_1)^2+(Y_1Z_2-Y_2Z_1)^2}{Z_1Z_2 \ \bar{||} \{ F_1, F_2 \}\bar{||} }(F_1\doth F_2),    \end{equation}
which  vanishes  if  $F_1\doth F_2=0$. \hfill $\square$


\subsection{Quadratic forms  on the symplectic plane}

Consider  the Poisson  algebra $\mathcal Q$  of the real  binary  quadratic  forms
\[   F=a p^2+2b pq+ c q^2, \]
provided with the  Poisson  bracket  in the  symplectic plane $(p,q)$.
We  will  rewrite  the form  $F$  as  linear  combination   of
\begin{equation}\label{base}    f_1:=2pq;  \quad  f_2:=p^2-q^2;  \quad   f_3:=p^2+q^2.   \end{equation}
From $F=Xf_1+Yf_2+Zf_3$  we obtain:
\begin{equation}\label{change}   X=b,  \quad  Y=(a-c)/2,  \quad  Z=(a+c)/2.         \end{equation}

Arnold  defined in \cite{Arn} a scalar  product $\tilde \Delta$,
as the symmetric  bilinear form  coinciding  on the diagonal
 with the discriminant $\Delta= ac-b^2$ of  the  form; if  $F=a p^2+2b pq+ c q^2$, and  $F'=a' p^2+2b' pq+ c' q^2$,  then
 \[  \tilde \Delta=  (ac'+ca')/2- bb'. \]

\noindent{\it Definition.}  We  define  the scalar product of two  forms
as   $-\tilde \Delta$.

\begin{prop} \label{prop2.5}  In the  new coordinates  the  scalar product of $F=(X,Y,Z)$ with $F'=(X',Y',Z')$
coincides  with  the  scalar product  (\ref{scalarph}) in  $\R^3_H$  and  the  square  norm  with  (\ref{normh}).
 \end{prop}

\noindent{\it Proof.} Using  \ref{change}, Eq. \ref{scalarph}  becomes
\[   F\doth F'= bb'-(ac'+ca')/2,  \]
 and, on the  diagonal:
   \[ \bar{||}F\bar{||} = b^2-ac.   \]
\hfill $\square$

\noindent{\it  Definition.} Two  forms  are said  {\it orthogonal}  iff  their  scalar product vanishes.

The  choice of the  base (\ref{base})  for  the quadratic  forms  allows  us  to write
the  Poisson  bracket  of  two  forms $F$ and $F'$  in  the   form analog to  the  vector  product.
\begin{prop}\label{prop2.6} The  algebra  $\Q$  is isomorphic to  the algebra  $\R^3_H$.
 \end{prop}
\noindent{\it  Proof.} In  the new  coordinates   the Poisson bracket of  $F$ and $F'$  is  four times
the  vector product (\ref{vectorph}) of  $F$ and  $F'$  in $\R^3_H$:
 \begin{equation}\label{eqPQ} \{F,F'\}=4 F\crossh F'.  \end{equation}
Indeed, the Poisson  brackets  between    the
basic functions  $f_i$,  $i=1,2,3$,  are  equal to:
\[  \{f_1,f_2\}=- \{f_2,f_1\}=-4 f_3, \quad    \{f_2,f_3\}= - \{f_3,f_2\}=4 f_1,\quad    \{f_3,f_1\}= - \{f_1,f_3\}=4 f_2.  \]
Eq. (\ref{eqPQ}) follows,  being  the  Poisson  bracket  bilinear  and antisymmetric.
The isomorphism $\mu: \R^3_H \rightarrow \Q$  associates  with the vector  $F$ the form  $\mu(F)=F/4$.   \hfill $\square$

\noindent{\it  Remark 9.}  The  functions  $f_1$, $f_2$  and  $f_3$  satisfy  $f_1^2+f_2^2-f_3^2=0$. Therefore
they can be  considered  as the coordinates in $\R^3$  of  a point  on the upper part ($Z\ge 0$) of the cone $X^2+Y^2-Z^2=0$.

\subsection{From  Algebra to Geometry}\label{hyperboloid}
Let us denote  by  $L^2_\pm$ the union of  the  two sheets  of  the two-sheeted hyperboloid.
The content  of   Section  \ref{sphere}   can  be  translated  for the hyperbolic
geometry  of $\R^3_H$.  In  fact,   the  definition of  the  duality  does  not depend  on the
   sign  of the  vector  product.     A  pair  of  antipodal  points on the  sphere  is  the intersection
   of a straight  line  through the origin  with the  sphere.   Here  it becomes a pair of antipodal points
either  in $L_\pm^2$ or in $D^2$,  whenever   the  straight line
does not  belong  to  the  cone.
The  great  circle  dual of a  pair  of  points is    the intersection of the plane orthogonal  (in  $\R^3_E$)
to the line  connecting  the  two  points.  Here   the great circle dual  of  a pair  of points    is replaced with  the  intersection
of the plane orthogonal (in $\R^3_H$) to the  line  connecting these  points  of the  hyperboloids.
Observe   that   such a plane intersects  either  only $D^2$  along  an ellipse  or  both $L_\pm^2$  and  $D^2$
along  a  pair  of hyperbolae.
The  conics obtained as  intersection of $L_\pm^2$  and  $D^2$ with the planes  through the origin   are  in fact the {\it  geodesics}  of  the  hyperboloids,  in the metrics (\ref{metricL})  and  (\ref{metricD}) on  $L^2$  and  $D^2$  respectively.

\begin{figure}[h]
\centerline{\epsfbox{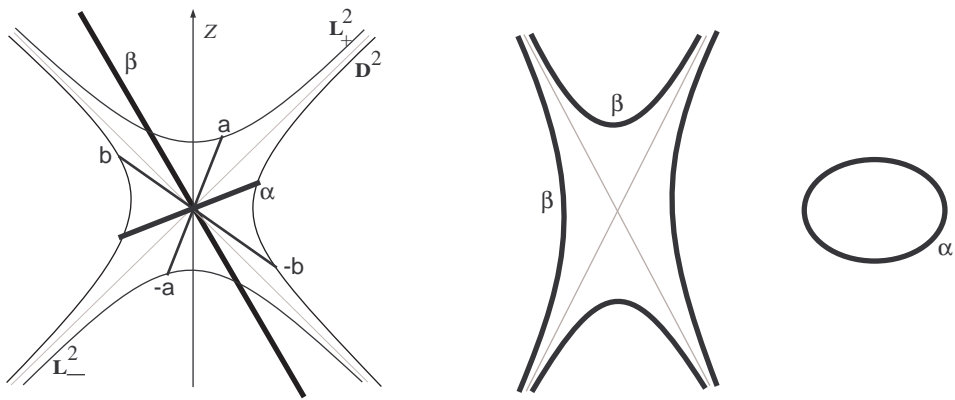}} \caption{Left:  section  $X=0$  of  $L_\pm^2\cup D^2$.  The  pair of   points
$(\at,-\at)$  and  $(\bt,-\bt)$ lie  on  the hyperboloids  and $\alpha$ and  $\beta$ are the geodesic  dual of them (Right).
   } \label{proj06}
\end{figure}

The  geodesic line $\alpha$ connecting  two points belonging to  $L_\pm^2\cup D^2$  is  either  an ellipse or a pair  of  hyperbolae (see  Figure  \ref{proj06}).  A  triangle is  a  subset of $L_\pm^2\cup D^2$
bounded by  three geodesics.

Also  here we  say  that  a triangle is  {\it proper}  if  no one of its  vertices  belong to  the pair of points dual  of  the  geodesics containing the other vertices,  or,  equivalently,
no  one of its  sides belongs to a geodesic dual  of  a pair containing the opposite vertex.

The  theorem  of the  altitudes  holds  here for proper  triangles    exactly  as  Proposition \ref{prop1.6}  and  Remark 6 in the spherical case,  substituting great circles with geodesics.

\section{Projective geometry}

We  consider now the projective plane $\R P^2$, with coordinates  $(x:= \frac{X}{Z},y:= \frac{Y}{Z})$,
to which the  space  $\R^3_H$  projects.

The  unit circle $x^2+y^2=1$ represents the  cone  $Z^2=X^2+Y^2$,  its interior  the Lobachevsky  disc  and the exterior
the De Sitter  world.

It is  easy  to  verify that  the  duality  between  a point  and  a line  (see  \cite{Arn})  on this projective
plane  is  the duality between   a vector based on  the origin and the  plane through the origin orthogonal in  $\R^3_H$ to
this vector.

\noindent{\it Remark 10.} From  Propositions \ref{prop2.1} and  \ref{prop2.3} it follows that the altitudes
theorem  obtained from  the  Jacobi  identity  of the Poisson  bracket for the  quadratic  forms algebra $\Q$
is  obtained  as well  from the Jacobi  identity  of  the  Poisson  bracket  for the  algebra
 $\H_1^+$ and for  the algebra  $\H^-_1$.

Recall  that   in $\R P^2$  we associate with   $F=(X,Y,Z)$,   the  pair  $(\at,\alpha)$, where  $\at=(X/Z,Y/Z)$  is the
projection of  $F$, and $\alpha$  is the line, projection  of the plane orthogonal  to  $F$.

\subsection{Relation  between $\R P^2_H$  and  $\R P^2_E$  }\label{relation}

Observe  that  the  orthogonality  between   two  lines, as well as the  duality between point  and line
in  the  projective  plane,  has not an intrinsic  meaning but is inherited by  the
orthogonality  in  $\R^3_H$  of their  preimages  by  the
projection.

In an analog way we can in fact  define the orthogonality  between   two  lines, as well as the  duality between point
and line   in the  projective  plane  as the
orthogonality in  $\R^3_E$ of their  preimages  by  the  projection.

\begin{figure}[h]
\centerline{\epsfbox{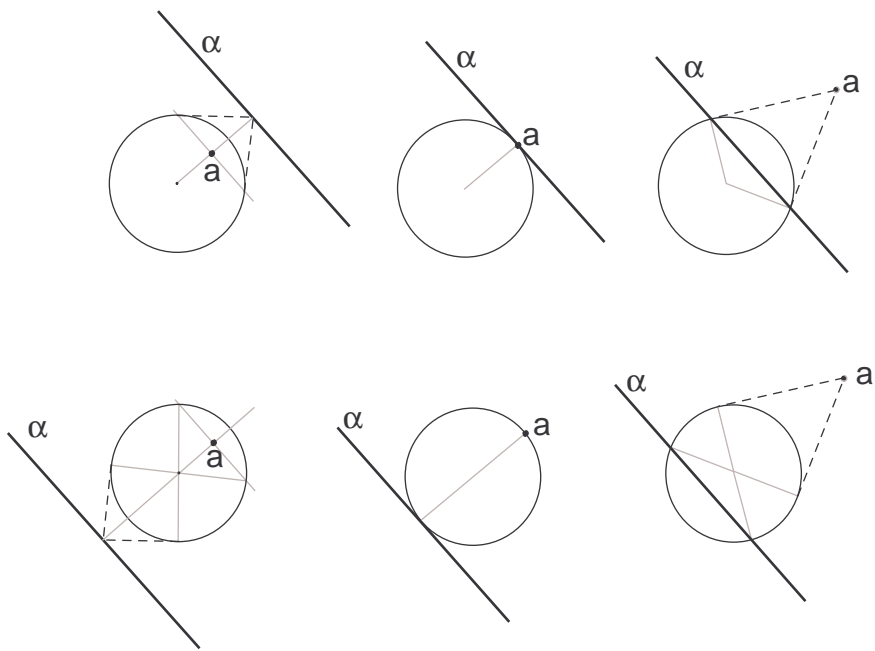}}  \caption{Three  different mutual  positions of a point $\at$ and  its dual line  $\alpha$. Upper  images:   duality   in $\R P^2_H$.  Lower images:  duality  in $\R P^2_E$.}\label{proj01}
\end{figure}

If we do  this, we  obtain  the  following  theorem, illustrated  in  Figure  \ref{proj01}.

Let  us denote  by  $\R P^2_H$  and  $\R P^2_E$  the  projective  plane  with the  orthogonality inherited  by
projectivization  from  $\R^3_H$ and  $\R^3_E$ respectively.  Provide both  planes  with
coordinates $(x:= \frac{X}{Z},y:= \frac{Y}{Z})$.

\begin{thm}\label{EH}  The  point dual  of a  line  in $\R P^2_E$  is  symmetric with respect  to the point $(0,0)$
of  the point dual  of the same line  in  $\R P^2_H$.  The line  dual  of  a point  in  $\R P^2_E$  is
symmetric  with respect to
the point $(0,0)$  to the  line  dual of the same point in  $\R P^2_H$.
\end{thm}

\begin{figure}[h]
\centerline{\epsfbox{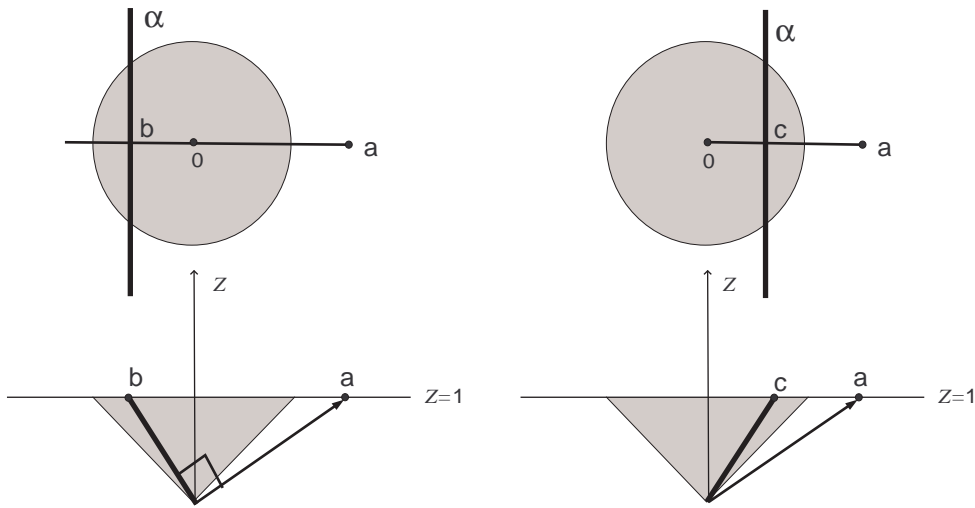}}  \caption{Left:   duality between the  point  $\at$ and  the line  $\alpha$  in $\R P^2_E$.  Right:  duality  between the point $\at$  and the line $\alpha$  in $\R P^2_H$.}\label{proj07}
\end{figure}

\noindent{\it Proof.}  Let  $\at=(a,0)$   and let  $\bt=(b,0)$ lie  on the  plane   orthogonal in  $\R^3_E$  to  $(a,0,1)$  (see  Fig. \ref{proj07}).     The vectors $(a,0,1)$  and $(b,0,1)$ satisfy  $ab+1=0$. Hence $b=-1/a$.  Let  $\ct=(c,0)$ lie  on the  plane   orthogonal in  $\R^3_H$  to  $(a,0,1)$.      The vectors $(a,0,1)$  and $(c,0,1)$ satisfy  $ac-1=0$. Hence $c=1/a$. Therefore  the points   $\bt$ and $\ct$  are symmetric  with respect to  the origin of $\R P^2$.

The  dual  line  $\alpha$ to $\at$  is  the line  through $\bt$   (in $\R P^2_E$) and  through  $\ct$
(in  $\R P^2_H$)  which is  orthogonal (in the plane  $Z=1$, with  Euclidean  metrics)
to  $\ot-\at$  (see Figure  \ref{proj07}).  \hfill  $\square$

\begin{cor} The Jacobi identity  for the Poisson  bracket  in the  algebra  $\S_1$
implies   the  concurrence  of the altitudes of the  triangle  in $\R P^2_E$.
\end{cor}

\begin{figure}[h]
\centerline{\epsfbox{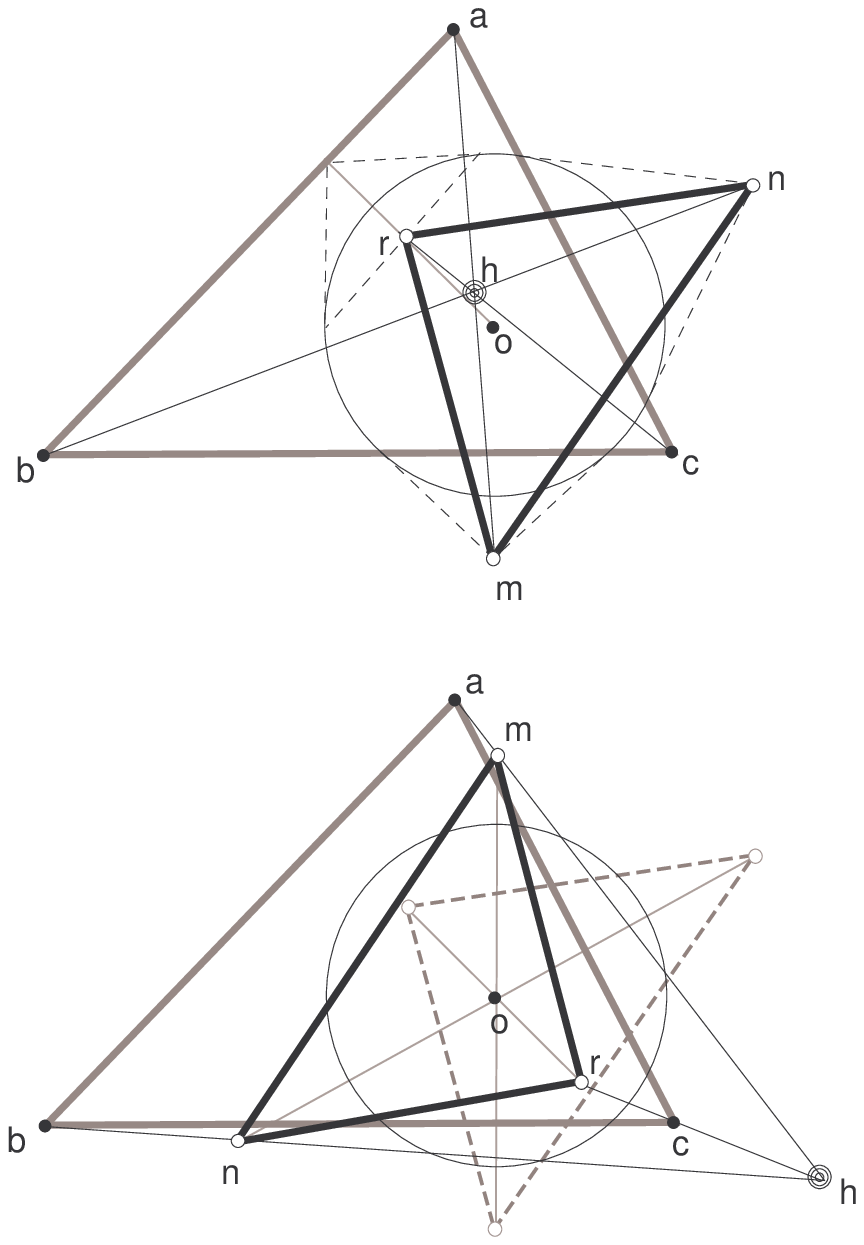}} \caption{ The triangle  $\at \bt \ct$  with its dual  triangle  $\mt \nt \rt$  in $\R P^2_H$
(upper image)  and  in  $\R P^2_E$ (lower image).  $\hh$  denotes  in both images   the  common point  of the  altitudes  of $\at \bt \ct$    and of   $\mt \nt \rt$.   The  points  $\mt, \nt, \rt$   in  $\R P^2_E$  are symmetric with respect  to the centre  $\ot$  of  the points    $\mt, \nt ,\rt$  in $\R P^2_H$.} \label{proj08}
\end{figure}

Figure  \ref{proj08}   shows the comparison of the altitudes  concurrence  theorem  for  a triangle in   $\R P^2_H$
and in $\R P^2_E$.

\subsection{Other models for  the projective  plane $\R P^2_H$ and  $\R P^2_E$  }
The projective  plane  $\R P^2_H$  can be  obtained  projecting  the  space  $\R^3_H$  on
any  affine plane  (not only  the plane $Z=1$).  This implies that  the  image of the  cone
 is  any conic section.

\begin{figure}[h]
\centerline{\epsfbox{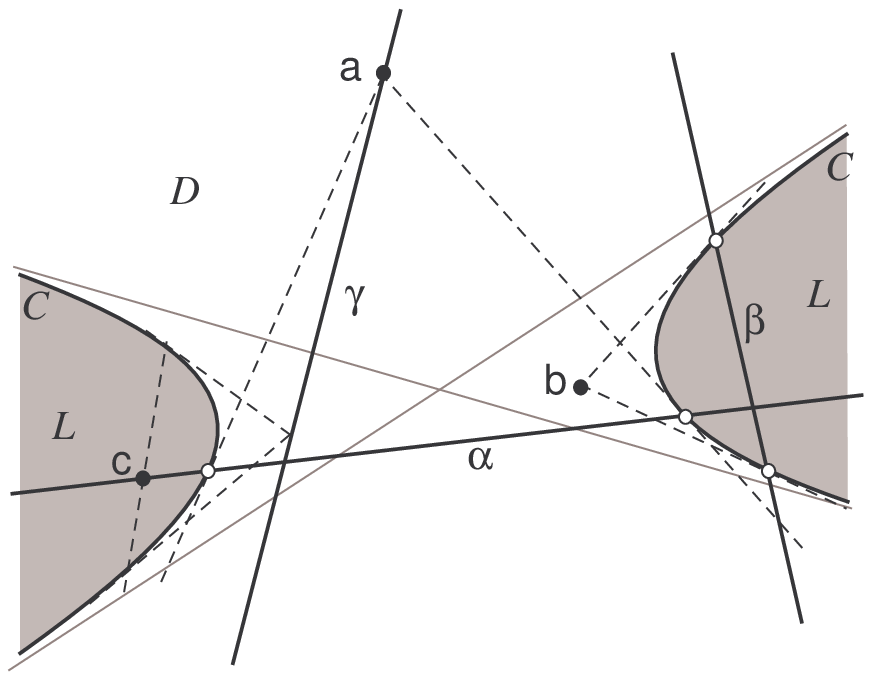}} \caption{The  lines $\alpha, \beta$ and  $\gamma$  are dual  of  the  points   $\at, \bt$ and  $\ct$ .  }\label {proj04}
\end{figure}

\noindent{\it  Definition.}   Given a  non degenerate  conic $C$  on the  projective plane,
define  the Lobachevsky  disc ($\L$) as   the  region of the complement  containing  the  foci  of the conic,  and  de Sitter
world ($\D$) the  other  region.  The  line  dual  of  a point $\at \in \D$  is  the line
joining the tangent  points  to  $C$  from $\at$.   The  point  dual  of  a line $\alpha$ intersecting  $C$
  is the   point of  $\D$   where  the  tangents  to $C$ at the points where $C$ intersects  $\alpha$ meet.      The  line  dual  of  a point $\at\in \L$ is  the line  containing all  points of $\D$ dual of   all lines  through
   $\at$ (see  Figure \ref{proj04}).  The conic $C$  is called  the {\it absolute}.

\noindent{\it  Definition.}  The  altitude  issued  from  a point
to a line  not containing it  is  the  line joining the point  with   the point,  dual
of the line.

\begin{prop}  Given any  non  degenerate conic  in a  plane, the  theorems of  projective   geometry    hold   in
that plane  taking the  above  definitions  of  duality.  \end{prop}

\begin{figure}[h]
\centerline{\epsfbox{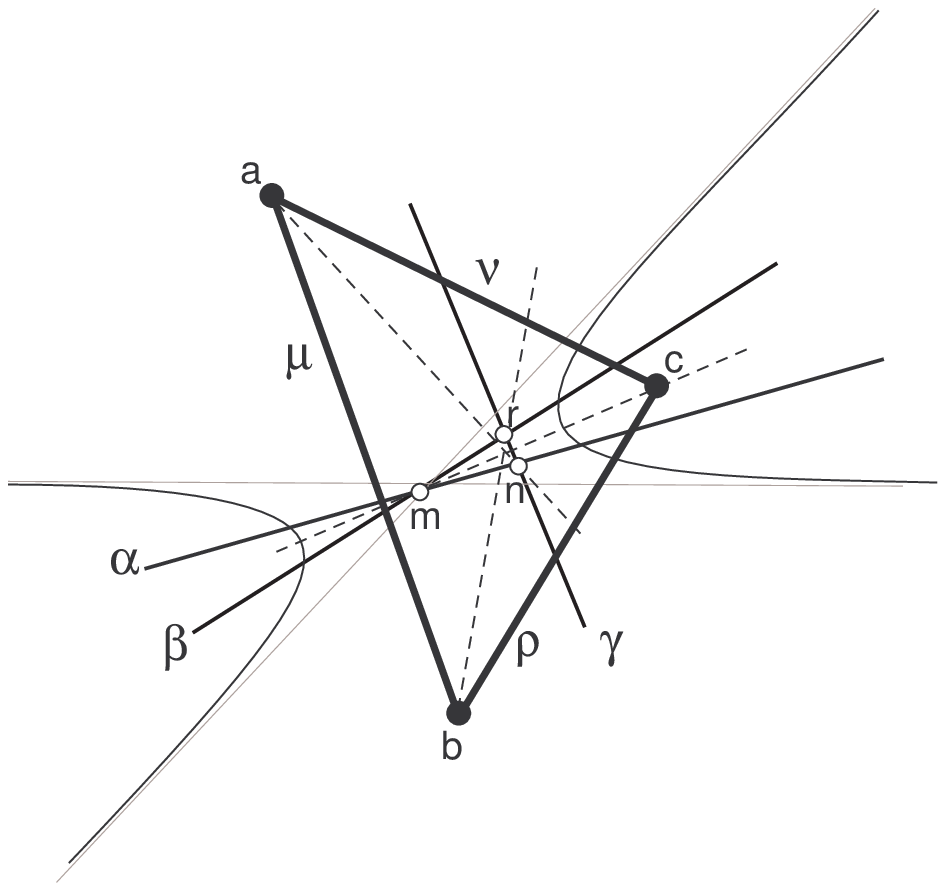}} \caption{The altitudes  of the  triangle  $\at \bt \ct$  meet  at the same  points.
They are as well  the altitudes of  the dual  triangle $\mt \nt \rt$.   } \label{proj05}
\end{figure}

\noindent{\bf  Example.}  An  illustration of the  hyperbolic  altitudes  concurrence theorem is   shown in  Figure    \ref{proj05}.

In  section \ref{relation}  we  have used the   constructions  of  dual  objects   in  $\R P^2_H$ to  obtain
  the dual  objects  in  $\R  P^2_E$.  In fact,   to  do  this  we  have  reflected points and  lines
  about  the origin of the  projective  plane.

  If   we  consider,  as a  model of  the  projective plane,  a  projection  of $\R^3$  to an affine plane,
  different  from  the plane $Z=1$,  the   point  $O$ about which we  have to  ``reflect''  points  and  lines, to obtain the  $\R P^2_E$ duality,
  becomes  the intersection  point  of that  plane with  the axis  of the  cone.  By  {\it reflection of a point $c$, lying
   on  $C$,}   we  mean  the second meeting point  of $C$ with
  $Oc$.  By   {\it reflection of a line   through two points} $c_1$ and $c_2$ of  the  conic $C$ with respect to  $O$   we  mean  the  line  connecting the
   points, different from $c_1$ and $c_2$, intersection of the conic $C$ with the  lines $Oc_1$  and $Oc_2$.
 By  {\it reflection of a point $p$, not lying  on  $C$,} with respect to $O$,   we  mean  the meeting point  of  two  lines,
  obtained  by reflection with respect to $O$  of two  lines meeting  at $p$.
 Finally,  by  {\it reflection of a  line $m$ nonintersecting  $C$} we  mean   the line  connecting two  points,
 which are  obtained by  reflection  of  two points  of   $m$.

      The  point $O$ is well  defined by the  conic itself  representing the absolute.  This is  however  an  indirect  and  complicated way  to define duality in  $\R P^2_E$.

 Indeed,  if  we  are  interested in the  geometry  of   $\R  P^2_E$,  obtained  as  projection  of $\R^3_E$ on  an
 arbitrary affine plane,  for the $\SO(3)$  invariance we can  take that  plane  as  the  plane  $Z=1$,
  the intersection  of this plane with the so  defined   $Z$  axis as the origin, and  the unit
  circle centered  at  the origin  as the  absolute,  so that Theorem \ref{EH}  holds.

\noindent{\it Remark 11.} The altitudes  concurrence theorem   is  mainly stated  for  triangles  in  the  Euclidean  plane.
In fact,  in this  case the  orthogonality   between  an altitude  and the relative  side  of   a  triangle  is
a  metric property of the Euclidean plane,  property that  is not invariant  under  projective transformations.  In   the  projective plane,  the  altitude  from  a vertex relative to  the opposite side of a triangle   is the straight line  connecting    that
vertex  with the  point,  dual  to  the  line  containing the opposite side. The   'orthogonality' between altitude and side  stands for  the  orthogonality (defined  as the vanishing of a  scalar product in  a  3-dimensional  metric space) of the corresponding planes  through the  origin   containing  respectively the  altitude  and the side, and is  therefore invariant  under  projective  transformations.

\section{Tomihisa's  identity}

The  following  identity  was  found  by  Tomihisa \cite{Tom} for  the  algebra  $\Q$.

Given  5  quadratic forms  $F_1, F_2, F_3, F_4, F_5$,
\begin{equation}\label{Tomi}
\{F_1,\{\{F_2,F_3 \},\{F_4, F_5 \}\}\}+\{F_3,\{\{F_2,F_5 \},\{F_4 , F_1 \}\}\}+\{F_5,\{ \{F_2,F_1 \},\{F_4 , F_3 \}\}\}=0.
\end{equation}

  The Tomihisa identity has  been  written  here   (\ref{Tomi})
in a little different form  with  respect  to   the expression  in \cite{Tom}   to  point out the    following  symmetry:
the  indices  (1,3,5)  are  cyclically  permuted  in the  three  terms
(like  the three indices of  the  Jacobi identity).  The  indices  (2,4)  are  at the same  place in all  terms.
The  indices  playing  this  role  will be  called  {\it  fixed  indices}
 of  the  identity.

We  have seen  that  the  algebra  $\Q$  is isomorphic to  the algebra $\R^3_H$.  We  will see  in  the last subsection   that   the Tomihisa identity holds  in all real  Lie-algebras
of  dimension 3,  and  hence  in all  Poisson  algebras  we  have  here  considered. It  has indeed  a unique  geometrical  meaning using  the  geometrical  interpretation  of  the Poisson  bracket
given in both sections  \ref{sphere} and  \ref{hyperboloid}.

\subsection{Geometrical  remarks  about  Tomihisa  identity}  An element of  the  Poisson algebras we have  considered   is  associated with  a pair  of  dual  geometric objects.
In this  way,  any  expression  concerning  Poisson  brackets  has  two  geometrical  meanings.

Moreover,  we   underline  the  following  fact. Call  $p(F)$  and  $\lambda(F)$ the  point  and its dual line
associated with  $F$.        The  expression:
\[   \{F,F'\}= F''  \]
has the following  meaning:
\begin{enumerate}
 \item $\lambda(F'')$ is the line  connecting  $p(F)$  and $p(F')$,
\item $p(F'')$  is  the intersection of  $\lambda(F)$  with  $\lambda(F')$,
\item $\lambda(F'')$  is  the  line  from  $p(F)$   orthogonal to $\lambda(F')$,
\item $\lambda(F'')$  is  the  line  from  $p(F')$  orthogonal to $\lambda(F)$.
\end{enumerate}

 We  say that  the  interpretations (1)  and (2) are {\it homogeneous}, since  the
objects inside the brackets  are  both  points or    both lines.

\noindent{\it Definition.} An expression  concerning $n$  objects  and  many levels of Poisson brackets
 is said  {\it homogeneous} if  there  exists a homogeneous interpretation
of  each   bracket,  whenever   the  $n$  objects are  interpreted  either all  as  points, either  all
as the  dual  lines.

\noindent{\it Example.}  The  Jacobi identity is  not  homogeneous,  since  if $F_1,F_2$ and $F_3$ are
interpreted  as points, $\{\{F_1,F_2\},F_3\}$ is  not  homogeneous.

{\it  Remark 12.} The  straight line  connecting two points and the intersection point of  two  lines
are  not  depending  on the orthogonality  definition.  Hence   a homogeneous  expression of the  considered  Lie  algebras  has  a  geometrical  meaning  independent of the definition  of duality  (i.e.,  of the definition of   orthogonality).

This  leads  us to the conclusion  that  the geometric  implications  of  the  Tomihisa  identity
must hold  for the projectivizations of  all  spaces  we have   considered,  independently  of their
metric structure.

\begin{figure}[h]
\centerline{\epsfbox{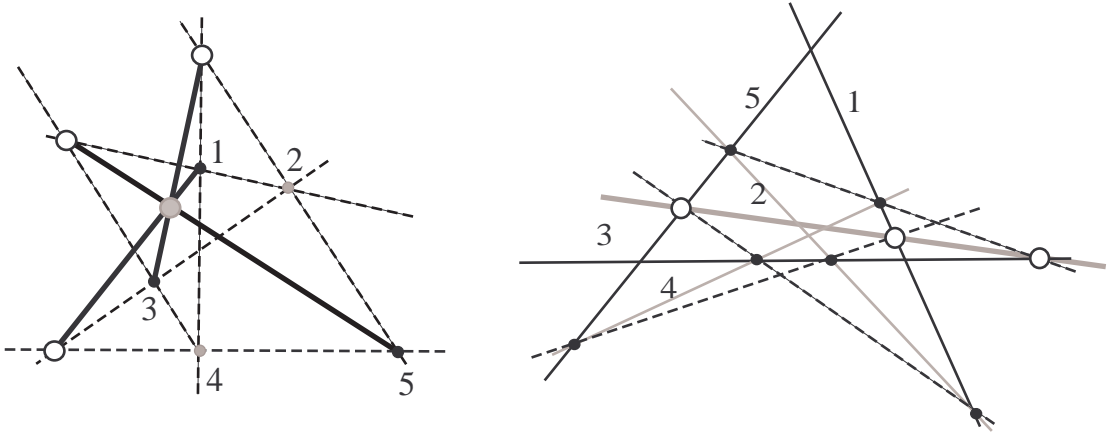}} \caption{ Tomihisa identity for  five elements  $F_i$.   The numbers 1-5 indicates
the indices $i$.  Left: the  $F_i$  are  interpreted as points: the  three bold lines intersect at the same point.
Right:  the  $F_i$ are  interpreted  as lines: the  three points marked by withe discs  lie  on the same
 line} \label{proj09}
\end{figure}

Figures  \ref{proj09} and \ref{proj10} illustrate the  Tomihisa identity  in $\R P^2$.

We  conclude  this  subsection with  a  comment  to   Remark  12  due to
V.Timorin.

\noindent {\it  Remark 13.} (V. Timorin) For a real oriented  vector space $V$ of dimension 3,
we can  define a  bilinear skew-symmetric map
$V\times V\to V^*$  (the dual space) in the  following  way.
Let    $\u,\v,\z$   be  three  vectors  of $V$.
Denote  by $P(\u,\v,\z)$  the oriented  volume spanned  by  the ordered  triple of   vectors $\u,\v,\z$.
The volume  has  a  sign which is  invariant under  cyclic  permutations  of the three vectors
and changes  sign under  odd  permutations  of them.
We define a linear functional $\mathcal P_{\u,\v}$  on
$V$  (i.e. an element of the dual space $V^*$) sending  $\z$  to  $P(\u,\v,\z)$.
The  map:  $V\times V\to V^*$,  sending  $(\u,\v)$  to   $\mathcal P_{\u,\v}$
is therefore  bilinear  and  skew-symmetric.
Since the volume form is well-defined up to a constant factor,
$\mathcal P$ is also well-defined up to a constant factor.
It is completely independent of any kind of additional structure,
say, orthogonality, duality between points and lines in the same
projective plane etc.
Similarly,   if we take two vectors $\x^*$ and $\y^*$  in the dual space,
then the  linear functional  $\mathcal P_{\x^*,\y^*}$  on  $V^*$, which sends $\z^*$
to  $P(\x^*,\y^*,\z^*)$,  is a vector in $V$ (since $V^{**}$ is
canonically isomorphic to $V$).
Thus  the Tomihisa identity  holds  interpreting  the  Poisson  bracket  as  the  functional $\mathcal P$,
without any reference  to  any Poisson structure.

\subsection{Tomihisa's  identity  as Pappus' theorem} In  \cite{Tom},  the Pappus theorem  is  deduced  from  the  identity (\ref{Tomi}) and other equations
\footnote{In particular,  in  \cite{Tom},  the Pappus Theorem is deduced  by an  identity (Theorem 2)   already  shown in \cite{Bla} as equivalent to the Pappus  Theorem},
 holding  for   the scalar product in the  space  of quadratic  forms.

 We  show here  that the identity  (\ref{Tomi}) contains  in fact  in itself the    Pappus  theorem.

\begin{prop} The Tomihisa identity for  the  points $F_1,\dots, F_5$  is equivalent  to the Pappus
theorem.
\end{prop}

{\it Proof.} Pappus' theorem states that  given  two  triples $(a_1,a_2,a_3)$  and $(b_1,b_2,b_3)$ of points,  each
one lying on a straight line, and defining   the point  $c_i$  ($i=1,2,3$)  as  the  meeting of the two straight  lines through $a_j,b_k$  and  $a_k,b_j$  ($i,j,k$ all different) then  the  points  $c_1,c_2,c_3$ are collinear.

\begin{figure}[h]
\centerline{\epsfbox{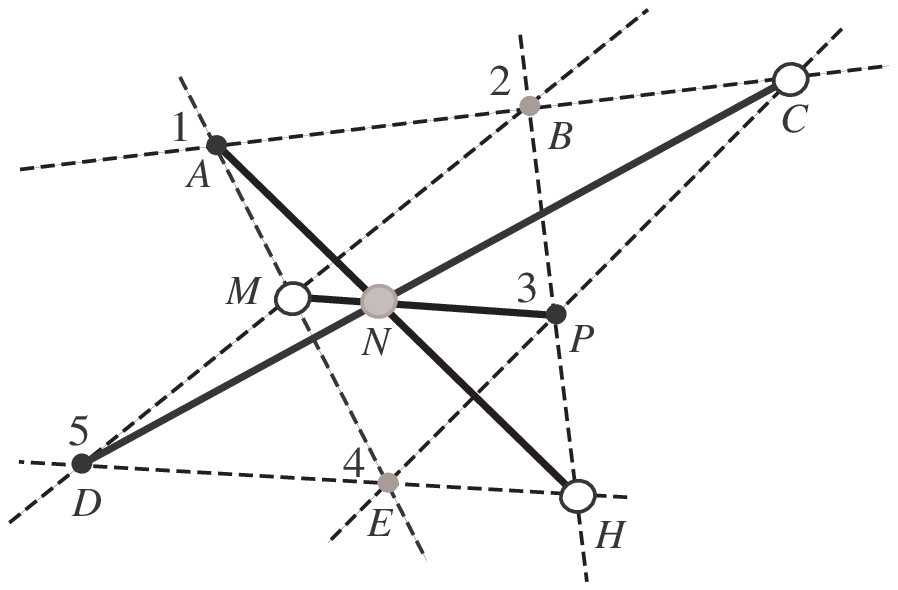}} \caption{}  \label{proj10}
\end{figure}

In  Figure  \ref{proj10} the  points  $a_i$  are named  $A,B,C$,  the points  $b_i$  are  named $D,E,H$ and the points
$c_i$  are  named  $M,N,P$.

An  equivalent  formulation  that the  points  $M,N,P$ are  collinear is   that  the  straight  line joining
$M$  and $P$  passes  through  $N$, the meeting point  of
the  lines connecting $A$  with $H$  and  $D$ with $C$.

Suppose  the points  $A,B,C$, $D,E,H$  be  assigned.    $P$  is the  meeting point  of $BH$  with $CE$.  We  set
\[ F_1\equiv A, \quad   F_2 \equiv B,  \quad   F_3\equiv P, \quad  F_4 \equiv E,  \quad F_5 \equiv D. \]
The  points  $C$  and  $H$  are therefore
\[  C= \{\{F_1,F_2\},\{F_3,F_4\}\},  \quad  H=\{\{F_3,F_2\},\{F_5,F_4\}\}.\]
Moreover
 \[ M=\{\{F_5,F_2\},\{F_1,F_4\}\}.\]
The  Tomihisa identity  is thus
\[   \{F_1,H\}+\{F_3,M\}+\{F_5,C\}=0,  \]
i.e.,   the straight  line connecting  $P$ and  $M$ passes  through the meeting point  of
the  lines connecting $A$  with $H$  and  $D$ with $C$.  \hfill $\square$

\subsection{An algebraic  remark about  Tomihisa's  identity}

According  to \cite{Raz}, there  are 2  identities  of  degree 5  that  constitute a  base  for the
identities  in  the simple real  Lie algebras  of dimension  3 (The vector space $\R^3$  and  $sl_2(\R))$):
Given  5  elements $X_0,X_1,X_2,X_3$ and $X_4$,  the following  identities hold:
\begin{equation}\label{ident1}
\sum_{s \in  S_4}  \sigma(s)[ X_{s(1)}, [ X_{s(2)}, [ X_{s(3)}, [ X_{s(4)}, X_0 ] ]    ] ]=0,
\end{equation}
where  $\sigma(s)$  is the  sign of  the  permutation $s$, positive  or  negative respectively   for  even or
odd  permutations;  and
\begin{equation}\label{ident2}
[ X_{0}, [ X_{0}, [ X_{0}, [ X_{1}, X_2 ] ]    ] ] - [ X_{1}, [ X_{0}, [ X_{0}, [ X_{2}, X_0 ] ]    ] ] + [ X_{2}, [ X_{0}, [ X_{0}, [ X_{0}, X_1 ] ]    ] ]=0. \end{equation}

\begin{prop} \label{Tomis}  Identity (\ref{ident1}) follows   from the  Tomihisa  identity.
\end{prop}

{\it  Proof.}  Write  Identity  (\ref{Tomi})  renaming  $F_2$  and
$F_4$,   respectively  $X_1$  and  $X_0$,  and  renaming  $F_1, F_3,
F_5$,  respectively,  $X_2,X_3,X_4$. The identity itself  in the  Lie  algebra will be
denoted  by

$T(X_0, X_1;  X_2,X_3,X_4)=0$.

Each  term  of  (\ref{Tomi}) is  of  type $[X_i,[[X_1,X_j],[X_0,X_k]]]$,  where  $(i,j,k)$  is
a  cyclic  permutation of  $(2,3,4)$.  We  apply  to  each  term
the  Jacobi  identity  in  this  way:
\[ [X_i,[[X_1,X_j],[X_0,X_k]]]=  [X_i, [[X_j,[X_k,X_0]],X_1]] +  [X_i,[[[X_k,X_0],X_1],
X_j]].\]  We   obtain:
\[
T(X_0, X_1;  X_2,  X_3, X_4)=
[X_2,[X_3,[X_1,[X_4,X_0]]]] -[X_2,[X_1,[X_3,[X_4,X_0]]]]+
\]
\[ [X_3,[X_4,[X_1,[X_2,X_0]]]]-[X_3,[X_1,[X_4,[X_2,X_0]]]]+
\]
\[+[X_4,[X_2,[X_1,[X_3,X_0]]]]-[X_4,[X_1,[X_2,[X_3,X_0]]]].
\]
Observe  that  in all  such terms  $X_0$  occupies  the  right position  of  the  fourth
Lie bracket,  in all  positive  terms  the  element  with  fixed index   $1$
occupies  the  left   position  in  the  third  Lie  bracket,   and  in  the  negative  terms
it  occupies  the  left   position of  the  second Lie  bracket.
The  other  elements  $X_2,X_3,X_4$  are  cyclically  permuted  in
 the  remaining  places,  so  that  the  sign of  the  term $ [X_i,[X_j,[X_k,[X_l,X_0]]]]$ coincides  with  the
  sign  of  the  permutation $(i,j,k,l)$  of  $(4,3,2,1)$.
The analogue terms  obtained  from  the Tomihisa identities  $T(X_0,X_i; X_j,X_k,X_l)$,
  with  fixed  indices $0,i$ for  $i=2,3,4$,
and  such  that $(i,j,k,l)$  is  an  even permutation of  $(4,3,2,1)$,  are therefore all  different.  We  obtain:

\[ \sum_{s \in  S_4}  \sigma(s)[ X_{s(1)}, [ X_{s(2)}, [ X_{s(3)}, [ X_{s(4)}, X_0 ] ]    ]
]=  T(X_0,X_1;  X_2,X_3,X_4)  \]
\[+
T(X_0, X_2 ;  X_3, X_1, X_4) +   T(X_0, X_3;  X_1, X_2, X_4)+
T(X_0, X_4;  X_2,  X_1, X_3).  \]

\hfill {$\square$}

{\it  Remark 14.}  In  fact,  A.  Dzhumadil'daev  has  recently  proved \cite{Dzhu} that   also  identity (\ref{ident2})  follows
from  Tomihisa  identity,  that therefore
constitutes itself   a  base.   He  provides  also in \cite{Dzhu} the  expression
of  the  Tomihisa  identity  using  the  base
$\{$(\ref{ident1}),(\ref{ident2})$\}$.

{\it  Remark 15.}  In  the vector
algebras  we  have  considered,  identity (\ref{ident2})  reads  as  the
Jacobi  identity  multiplied  by  the  square  norm  of  $X_0$. Its    geometrical  meaning
coincides therefore  with  that  of  the  Jacobi  identity.

Acknowledgements. I  express  my  gratitude  to V.I.Arnold,  who  suggested  to search  for the relation between
Tomihisa  identity and known  identities  in Lie algebras.   I am  also  grateful  to the  Referee, to  V. Timorin, to A.  Dzhumadil'daev
and to E. Malkovich  for their interest and remarks.

\end{document}